# On Tate-Shafarevich Groups of some Elliptic Curves

*Franz Lemmermeyer*

**Abstract.** Generalizing results of Stroeker and Top we show that the 2-ranks of the Tate-Shafarevich groups of the elliptic curves $y^2 = (x+k)(x^2+k^2)$ can become arbitrarily large. We also present a conjecture on the rank of the Selmer groups attached to rational 2-isogenies of elliptic curves.

## 1. Introduction

In order to introduce the relevant notation, we start by reviewing the well known first descent on elliptic curves by rational 2-isogenies (see e.g. [2, 7, 14]). Such isogenies of an elliptic curve $E$ defined over $\mathbb{Q}$ exist if and only if $E$ has a rational point of order 2, and by a suitable choice of the coordinate system we may assume that $E$ has a Weierstraß model of the form $E : y^2 = x(x^2 + ax + b)$ with $a, b \in \mathbb{Z}$. The presence of a torsion point $T = (0,0)$ of order 2 guarantees the existence of a rational 2-isogeny $\phi : E \longrightarrow \widehat{E}$ onto $\widehat{E} : y^2 = x(x^2 + \widehat{a}\,x + \widehat{b})$, where $\widehat{a} = -2a$ and $\widehat{b} = a^2 - 4b$. There is a dual isogeny $\psi : \widehat{E} \longrightarrow E$ such that $\psi \circ \phi = [2]_E$, where $[2]_E$ denotes the multiplication by 2 on $E$; similarly, $\phi \circ \psi$ is multiplication by 2 on $\widehat{E}$.

Weil [16] studied the map $\alpha$ defined by $\alpha(P) = x\mathbb{Q}^{\times 2}$ for all points $P = (x,y) \in E(\mathbb{Q})$ different from $\mathcal{O}$ and $T$, and showed that $\alpha$ becomes a homomorphism $E(\mathbb{Q}) \longrightarrow \mathbb{Q}^\times / \mathbb{Q}^{\times 2}$ by putting $\alpha(\mathcal{O}) = \mathbb{Q}^{\times 2}$ and $\alpha(T) = b\mathbb{Q}^{\times 2}$. The kernel of $\alpha$ coincides with the image of the 2-isogeny $\psi$ introduced above. In other words, we have the exact sequence

$$0 \longrightarrow \psi(\widehat{E}(\mathbb{Q})) \longrightarrow E(\mathbb{Q}) \xrightarrow{\alpha} \mathbb{Q}^\times/\mathbb{Q}^{\times 2}$$

as well as the corresponding result for the dual isogeny:

$$0 \longrightarrow \phi(E(\mathbb{Q})) \longrightarrow \widehat{E}(\mathbb{Q}) \xrightarrow{\beta} \mathbb{Q}^\times/\mathbb{Q}^{\times 2}.$$

Moreover Tate showed that the $\mathbb{Z}$-rank $r$ of the Mordell-Weil groups $E(\mathbb{Q})$ and $\widehat{E}(\mathbb{Q})$ is given by

$$2^{r+2} = \#\operatorname{im}\alpha \cdot \#\operatorname{im}\beta. \tag{1}$$

The images of $\alpha$ and $\beta$ can be described explicitly as follows: $W(\widehat{E}/\mathbb{Q}) := \operatorname{im}\alpha$ (note that $\operatorname{im}\alpha \simeq E(\mathbb{Q})/\psi(\widehat{E}(\mathbb{Q}))$) consists of all classes $b_1 \mathbb{Q}^{\times 2}$, where $b_1$ is a squarefree integer with $b_1 b_2 = b$, such that

$$N^2 = b_1 M^4 + aM^2 e^2 + b_2 e^4 \tag{2}$$

has a nontrivial[1] primitive[2] solution in integers $N, M, e \in \mathbb{N}$. The equation (2) is called a *torsor* of $E/\mathbb{Q}$ and will be denoted by $\mathcal{T}^{(\psi)}(b_1)$; note that these torsors depend on the 2-isogeny $\psi : \widehat{E} \longrightarrow E$: elliptic curves with three rational points of order 2 have three 2-isogenies. It is easy to see that every rational point $P \neq \mathcal{O}$ on $E$ has the form $P = (m/e^2, n/e^3)$ for integers $n, m, e \in \mathbb{Z}$ such that $(m, e) = (n, e) = 1$, and we have $\alpha(P) = m\mathbb{Q}^{\times 2}$ by definition; moreover, it can be shown that the corresponding torsor $\mathcal{T}^{(\psi)}(m)$ is solvable. Conversely, if $(N, M, e)$ is a nontrivial primitive solution of $\mathcal{T}^{(\phi)}(b_1)$, then $(b_1 M^2/e^2, b_1 MN/e^3)$ is a rational point on $E$, and the group structure on the torsors of $W(\widehat{E}/\mathbb{Q})$ is induced by the addition law on the elliptic curve: pick two rational points on $E$ corresponding to rational solutions on the torsors $\mathcal{T}^{(\psi)}(\beta_1)$ and $\mathcal{T}^{(\psi)}(\beta_2)$, add the points on $E$, and map this sum to the corresponding torsor; it turns out that this torsor is $\mathcal{T}^{(\psi)}(\beta_3)$, $\beta_3$ being the square free kernel of $\beta_1 \beta_2$; thus the multiplication on the set of torsors coincides with the multiplication in $\mathbb{Q}^{\times}/\mathbb{Q}^{\times 2}$. We also remark that rational solutions of $\mathcal{T}^{(\psi)}(b_1)$ can be made integral and primitive by clearing denominators etc.

Similarly, $W(E/\mathbb{Q}) := \operatorname{im}\beta$ consists of all classes $b_1 \mathbb{Q}^{\times 2}$, $b_1$ squarefree, with $b_1 b_2 = \widehat{b} = a^2 - 4b$ such that the torsor

$$\mathcal{T}^{(\phi)}(b_1) : N^2 = b_1 M^4 - 2aM^2 e^2 + b_2 e^4 \tag{3}$$

has a nontrivial primitive integral solution.

The classes $b_1 \mathbb{Q}^{\times 2}$ such that (3) has solutions in every $\mathbb{Q}_p$ (including $\mathbb{R} = \mathbb{Q}_\infty$) form a subgroup $S^{(\phi)}(E/\mathbb{Q})$ of $\mathbb{Q}^{\times}/\mathbb{Q}^{\times 2}$ (called the $\phi$-part of Selmer group of $E$) that contains $W(E/\mathbb{Q})$ as a subgroup. The corresponding factor group is called a Tate-Shafarevich group; thus we have, by definition, the exact sequences

$$\begin{array}{ccccccccc}
0 & \longrightarrow & W(E/\mathbb{Q}) & \longrightarrow & S^{(\phi)}(E/\mathbb{Q}) & \longrightarrow & \mathrm{III}(E/\mathbb{Q})[\phi] & \longrightarrow & 0, \\
0 & \longrightarrow & W(\widehat{E}/\mathbb{Q}) & \longrightarrow & S^{(\psi)}(\widehat{E}/\mathbb{Q}) & \longrightarrow & \mathrm{III}(\widehat{E}/\mathbb{Q})[\psi] & \longrightarrow & 0.
\end{array} \tag{4}$$

Selmer and Tate-Shafarevich groups can be defined for any isogeny of elliptic curves; moreover it can be shown that $\mathrm{III}(E/\mathbb{Q})[\phi]$ injects into $\mathrm{III}(E/\mathbb{Q})[2]$ (see the diagram in Section 5 below). It is known that there are Tate-Shafarevich groups with arbitrarily large cardinality: Cassels [4] used the pairing on $\mathrm{III}(E/\mathbb{Q})[3]$ for showing that there are elliptic curves $E$ defined over $\mathbb{Q}$ such that the 3-rank of $\mathrm{III}(E/\mathbb{Q})[3]$ is arbitrarily large; a similar method was used by McGuinness [12] to prove the same for $\mathrm{III}(E/\mathbb{Q})[2]$. Bölling [3] and Kramer [8] also showed that $\mathrm{III}(E/\mathbb{Q})[2]$ can be arbitrarily large. In all these constructions, reciprocity plays an important role; on the other hand, these authors also make use of heavier machinery such as the Cassels pairing on $\mathrm{III}$. In this paper we want to show that large $\mathrm{III}$ can be found with reciprocity alone. This is accomplished by factoring the torsors (2) and (3) over the quadratic number fields $\mathbb{Q}(\widehat{E}[2])$ and $\mathbb{Q}(E[2])$, respectively: multiplying the torsor $\mathcal{T}^{(\psi)}(b_1)$ in (2) by $b_1$ we get, for example,

$$b_1 N^2 = (b_1 M + \tfrac{1}{2}ae^2)^2 - \tfrac{1}{4}(a^2 - 4b)e^4, \qquad (5)$$

and the right hand side splits over $\mathbb{Q}(\sqrt{a^2 - 4b}) = \mathbb{Q}(E[2])$. The idea is then to use the arithmetic of $\mathbb{Q}(E[2])$ to deduce necessary conditions for the existence of a rational point on $\mathcal{T}^{(\psi)}(b_1)$.

---

[1] That is, $N = M = e = 0$ is excluded.
[2] This is our abbreviation for $(N, e) = (M, e) = 1$.

**Remark.** The torsor (5) splits not only over $\mathbb{Q}(\widehat{E}[2])$; writing it in the form $\tfrac{1}{4}(a^2 - 4b)e^4 = (b_1 M + \tfrac{1}{2}ae^2)^2 - b_1 N^2$ shows that it also splits over $\mathbb{Q}(\sqrt{b_1})$.

In this paper, we study the elliptic curve $y^2 = (x+k)(x^2+k^2)$ for integers $k \neq 0$. In fact, we will work with the model $E_k : y^2 = x(x^2 - 2kx + 2k^2)$ throughout this paper. The curve $\widehat{E}_k : y^2 = x(x^2 + 4kx - 4k^2)$ is the 2-isogenous curve of $E$. Both $E_k$ and $\widehat{E}_k$ have conductor $N = 2^7 k^2$, and their torsion groups have order 2 (see [13], but observe that one still has to check that the 3-torsion is trivial), the points of order 2 being $T = (0,0)$ and $\widehat{T} = (0,0)$, respectively. Note that $\alpha(T) = 2\mathbb{Q}^{\times 2}$ and $\beta(\widehat{T}) = -\mathbb{Q}^{\times 2}$, so $W(\widehat{E}_k/\mathbb{Q})$ and $W(E_k/\mathbb{Q})$ always have even cardinality. Also the fact that the fields $\mathbb{Q}(E_k[2]) = \mathbb{Q}(i)$ and $\mathbb{Q}(\widehat{E}_k[2]) = \mathbb{Q}(\sqrt{2})$ have class number 1 keeps things simple. Finally we observe that the same technique can be used to construct families of non-congruent numbers (see [11]).

## 2. Quadratic Residue Symbols

Let $\mathcal{O}_k$ be the ring of integers of a number field $k$; for prime ideals $\mathfrak{p}$ of $\mathcal{O}_k$ and integers $\alpha \in \mathcal{O}_k \setminus \mathfrak{p}$ with odd norm, define the quadratic residue symbol $[\alpha/\mathfrak{p}] = \pm 1$ by demanding that $[\alpha/\mathfrak{p}] \equiv \alpha^{(N\mathfrak{p}-1)/2} \bmod \mathfrak{p}$. For $k = \mathbb{Q}$, this is the usual Legendre symbol $(\cdot/\cdot)$. For quadratic number fields it follows directly from the definition that $[\alpha/p] = (N\alpha/p)$ for any rational prime $p$, and that $[a/\mathfrak{p}] = (a/p)$ whenever $a \in \mathbb{Z}$, where $\mathfrak{p}$ is a prime ideal in $k$ of degree 1 with norm $p$. For quadratic residues $a \equiv x^2 \bmod p$ of primes $p \equiv 1 \bmod 4$, the biquadratic symbol $(a/p)_4$ is defined by $(a/p)_4 = (x/p)$.

**Lemma 1.** *Let $k$ be a squarefree product of primes $p_j \equiv 1 \bmod 8$ with $(p_i/p_j) = +1$ for $i \neq j$, and let $\pi_j = e + f\sqrt{2}$ be an element of norm $p_j$ in $\mathbb{Z}[\sqrt{2}]$. If $k_1$ is a product of $p_j$ and $k_1 = \kappa_1 \overline{\kappa}_1$ for some $\kappa_1 \in \mathbb{Z}[\sqrt{2}]$, then $[\kappa_1/\pi]$ does not depend on the choice of the $\pi_j$, $1 \leq j \leq r$.*

*Proof.* It is clearly sufficient to show that $[\lambda/\pi] = [\overline{\lambda}/\pi]$, where $\lambda = \pi_i$ and $\pi = \pi_j$ with $i \neq j$, and where $\overline{\lambda}$ is the conjugate of $\lambda$. But this follows immediately from $[\lambda\overline{\lambda}/\pi] = [l/\pi] = (l/p) = +1$. □

**Lemma 2.** *Let* $k = \mathbb{Q}(\sqrt{-1})$ *and let* $\sigma$ *be a prime in* $\mathcal{O}_k$ *with norm* $\sigma\overline{\sigma} = p \equiv 1 \bmod 8$. *Then* $[\sigma/\overline{\sigma}] = 1$.

*Proof.* This is a consequence of the quadratic reciprocity law in $\mathbb{Z}[i]$ due to Gauß and Dirichlet (note that we may assume $\sigma \equiv 1 \bmod 2\mathcal{O}_k$ without loss of generality since $[i/\sigma] = (-1)^{(p-1)/4} = 1$). A direct proof goes as follows: assume that $m = -1$ and write $\sigma = a + bi$. Then $[\sigma/\overline{\sigma}] = [2a/\overline{\sigma}] = (2a/p) = 1$ since $(a/p) = (p/a) = 1$. □

**Lemma 3.** *Let* $k = \mathbb{Q}(\sqrt{2})$ *and let* $\pi \equiv 1 \bmod 2\mathcal{O}_k$ *be a prime in* $\mathcal{O}_k$ *with norm* $\pi\overline{\pi} = p \equiv 1 \bmod 8$. *Then* $[\pi/\overline{\pi}] = (2/p)_4$.

*Proof.* Write $\pi = e + f\sqrt{2}$. Then $[\pi/\overline{\pi}] = (2e/p) = (e/p)$. On the other hand, $(2/p)_4 = [\sqrt{2}/\pi]$, and since $\sqrt{2} \equiv -e/f \bmod \pi$, we have $[\sqrt{2}/\pi] = [-ef/\pi] = (-ef/p) = (e/p)$ since $(f/p) = 1$: in fact, write $f = 2^j f'$ for some odd $f'$; then $j \geq 2$ and $(f/p) = (f'/p) = (p/f') = (e/f')^2 = 1$. □

Incidentally, the same result holds in $\mathbb{Z}[\sqrt{-2}]$.

**Lemma 4.** *Let* $p, q \equiv 1 \bmod 8$ *be primes; then we can choose* $\pi, \lambda \in \mathbb{Z}[\sqrt{2}]$ *such that* $N\pi = p$, $N\lambda = q$ *and* $\pi \equiv \lambda \equiv 1 \bmod 2$, *and we have* $[\pi/\lambda] = [\lambda/\pi]$.

*Proof.* This is well known and follows from the general quadratic reciprocity law in number fields with odd class number proved by Hilbert, Hecke, Dörrie, Hasse and others (see [9]). In order to convince the reader that its proof is completely elementary, we will give it here.

Write $\pi = a + b\sqrt{2}$, $\lambda = c + d\sqrt{2}$, and $d = 2^j d'$ for some odd $d'$. Then $[d/\lambda] = (d/q) = (d'/q) = (q/d') = (c^2 - 2d^2/d') = (c/d')^2 = +1$. This implies that $[\pi/\lambda] = [d/\lambda][ad + bd\sqrt{2}/\lambda] = [ad - bc/\lambda] = (ad - bc/q)$ since $d\sqrt{2} \equiv -c \bmod \lambda$. Similarly, we get $[\lambda/\pi] = (ad - bc/p)$, hence it is sufficient to show that $(ad - bc/pq) = 1$. To this end, write $ad - bc = 2^j e$ and notice that $pq = (ac - 2bd)^2 - 2(ad - bc)^2 \equiv (ac - 2bd)^2 \bmod e$; this shows that $(ad - bc/pq) = (e/pq) = (pq/e) = (ac - 2bd/e)^2 = +1$ as desired. □

We also note that, for primes splitting in $\mathbb{Q}(\zeta_8)$, i.e. primes $p \equiv 1 \bmod 8$, the relation $1 + i = \sqrt{2}\zeta_8$ implies that $(\frac{1+i}{p}) = (2/p)_4(-1)^{(p-1)/8}$; in particular we have $(\frac{1+i}{p}) = -(2/p)_4$ if $p \equiv 9 \bmod 16$, and $(2/p)_4 = (-1)^{(p-1)/8}$ if $(\frac{1+i}{p}) = +1$. Here we have written $(\frac{1+i}{p})$ instead of $[\frac{1+i}{\sigma}]$, where $\sigma \in \mathbb{Z}[i]$ is an element with norm $p$; this is no problem since $[\frac{1+i}{\sigma}]$ does not depend on the choice of $\sigma$.

**Lemma 5.** *Let* $p \equiv 1 \bmod 8$ *be prime; then there exist integers* $a, b, e, f \in \mathbb{N}$, $b$ *and* $f$ *even, such that* $p = a^2 + b^2 = e^2 - 2f^2$. *Then* $(\frac{1+i}{p}) = (\frac{1+\sqrt{2}}{p}) = (\frac{2}{a+b}) = (-1)^{e+f}$.

*Proof.* The first equality follows from $(1 + \zeta_8)^2 = (1 + i)(1 + \sqrt{2})$, the equality $(\frac{1+i}{p}) = (\frac{2}{a+b})$ is well known and easy to prove, and finally we have

$$\Big(\frac{1+\sqrt{2}}{p}\Big) = \Big[\frac{1+\sqrt{2}}{e+f\sqrt{2}}\Big] = \Big[\frac{f}{e+f\sqrt{2}}\Big]\Big[\frac{f+f\sqrt{2}}{e+f\sqrt{2}}\Big] = \Big[\frac{f}{e+f\sqrt{2}}\Big]\Big[\frac{f-e}{e+f\sqrt{2}}\Big].$$

But $[f/e + f\sqrt{2}] = (f/p) = 1$: in fact, $f = 2^j f'$ for some odd integer $f'$ and $(f/p) = (f'/p) = (p/f') = 1$. Similarly, $[f - e/e + f\sqrt{2}] = (f - e/p) = (e - f/p)$ since $(-1/p) + 1$. Thus

$$\Big(\frac{1+\sqrt{2}}{p}\Big) = \Big(\frac{e-f}{p}\Big) = \Big(\frac{p}{e-f}\Big) = \Big(\frac{-1}{e+f}\Big),$$

where we have used that a) $p = e^2 - 2f^2 \equiv -f^2 \bmod e - f$ and b) $e - f > 0$ and $e - f \equiv e + f \bmod 4$. □

## 3. The Tate-Shafarevich group

In this section we study the curve $E_k : y^2 = x(x^2 - 2kx + 2k^2)$, and in particular the torsors

$$\begin{aligned} \mathcal{T}^{(\phi)}(b_1) : N^2 &= b_1 M^4 + 4kM^2 e^2 - 4k^2 e^4, \\ \mathcal{T}^{(\psi)}(b_1) : N^2 &= b_1 M^4 - 2kM^2 e^2 + 2k^2 e^4. \end{aligned}$$

The following Proposition describing the $\phi$- and $\psi$-part of the Selmer groups of $E_k$ and $\widehat{E}_k$ is a special case of a result proved by S. Schmitt in [13]:

**Proposition 6.** *Suppose that $k = p_1 \cdots p_t$ is a product of primes $p_i \equiv 1 \bmod 8$ such that $(p_i/p_j) = +1$ for $i \neq j$ and $(\frac{1+i}{p_j}) = +1$ for all $1 \leq j \leq t$; then*

$$\begin{aligned} S^{(\phi)}(E_k/\mathbb{Q}) &= \{b_1 \mathbb{Q}^{\times 2} : b_1 \text{ squarefree, } b_1 \mid 2k\}, \\ S^{(\psi)}(\widehat{E}_k/\mathbb{Q}) &= \{b_1 \mathbb{Q}^{\times 2} : b_1 \text{ squarefree, } 0 < b_1 \mid 2k\}. \end{aligned}$$

The case $t = 1$ was discussed by R. Stroeker & J. Top [15]. Similar results for products of primes not necessarily $\equiv 1 \bmod 8$ can be found in S. Schmitt [13].

**Proposition 7.** *Consider the curve $E = E_k$, where $k = p_1 \cdots p_t$ is a product of primes $p_i \equiv 1 \bmod 8$, and assume that $k \equiv 9 \bmod 16$ if $t$ is odd. Suppose moreover that $(p_i/p_j) = +1$ for $i \neq j$ and $(\frac{1+i}{p}) = +1$ for all $p \mid k$; write $p_i = \pi_i \overline{\pi}_i$ for elements $\pi_i, \overline{\pi}_i \in \mathbb{Z}[\sqrt{2}]$ with $\pi_i \equiv 1 \bmod 2$, and assume in addition that $[\pi_i/\pi_j] = -1$ for $i \neq j$. Then $W(E/\mathbb{Q}) = \{\mathbb{Q}^{\times 2}, -\mathbb{Q}^{\times 2}\}$ if $t$ is odd, and $W(E/\mathbb{Q}) = \{\mathbb{Q}^{\times 2}, -\mathbb{Q}^{\times 2}\}$ or $W(E/\mathbb{Q}) = \{\pm\mathbb{Q}^{\times 2}, \pm 2k\mathbb{Q}^{\times 2}\}$ if $t$ is even, the first possibility being prohibited by the parity conjecture.*

*In particular, we have $\#\mathrm{III}(E/\mathbb{Q})[\phi] = \frac{1}{2}\#S^{(\phi)}(E/\mathbb{Q}) = 2^{t+1}$ if $t$ is odd, and $\#\mathrm{III}(E/\mathbb{Q})[\phi] \geq \frac{1}{4}\#S^{(\phi)}(E/\mathbb{Q}) = 2^t$ if $t$ is even, with equality if the parity conjecture is true.*

A similar result holds for the isogenous curve:

**Proposition 8.** *Consider the curve $E = E_k$, where $k = p_1 \cdots p_t$ is a product of primes $p_i \equiv 1 \bmod 8$ with $(p_i/p_j) = +1$ for $i \neq j$ and $(\frac{1+i}{p}) = +1$ for all $p \mid k$; write $p_i = \sigma_i \overline{\sigma}_i$ for elements $\sigma_i, \overline{\sigma}_i \in \mathbb{Z}[i]$, and assume in addition that $[\sigma_i/\sigma_j] = -1$ for $i \neq j$ (these symbols do not depend on the choice of the $\sigma_i$). Then $W(\widehat{E}/\mathbb{Q}) = \{\mathbb{Q}^{\times 2}, 2\mathbb{Q}^{\times 2}\}$ if $t$ is even, and $W(\widehat{E}/\mathbb{Q}) \subseteq \{\mathbb{Q}^{\times 2}, 2\mathbb{Q}^{\times 2}, k\mathbb{Q}^{\times 2}, 2k\mathbb{Q}^{\times 2}\}$ if $t$ is odd (with equality if the parity conjecture is true).*

*In particular, we have $\#\mathrm{III}(\widehat{E}/\mathbb{Q})[\psi] \geq \frac{1}{4}\#S^{(\psi)}(\widehat{E}/\mathbb{Q}) = 2^{t-1}$ if $t$ is odd, and $\#\mathrm{III}(\widehat{E}/\mathbb{Q})[\psi] = \frac{1}{2}\#S^{(\psi)}(\widehat{E}/\mathbb{Q}) = 2^t$ if $t$ is even, with equality if the parity conjecture holds.*

Combining these two propositions we find that

**Theorem 9.** *If the conditions in Propositions 7 and 8 are satisfied, then $E_k(\mathbb{Q})$ and $\widehat{E}_k(\mathbb{Q})$ have rank $\leq 1$, with equality if the parity conjecture holds. In this case, we also have*

$$\#\mathrm{III}(E/\mathbb{Q})[\phi] = \begin{cases} 2^{t+1} & \text{if } t \text{ is odd,} \\ 2^t & \text{if } t \text{ is even,} \end{cases} \quad \text{and} \quad \#\mathrm{III}(\widehat{E}/\mathbb{Q})[\psi] = \begin{cases} 2^{t-1} & \text{if } t \text{ is odd,} \\ 2^t & \text{if } t \text{ is even.} \end{cases}$$

*Proof.* The inequality rank $E_k(\mathbb{Q}) \leq 1$ follows by applying Tate's formula (1). By results in [13], the parity conjecture implies that the curves $E_k$ have odd rank if $0 < k \equiv 1 \bmod 8$; since we have proved that the rank is at most 1, we must have rank $E_k(\mathbb{Q}) = 1$. This in turn implies that the Tate-Shafarevich groups in Propositions 7 and 8 have square order. □

**Remark.** Note that standard theorems in analytic number theory imply that the conditions for Propositions 7 and 8 can be satisfied simultaneously for infinitely many primes: just look at the splitting behaviour of primes in the multiquadratic extensions of $k = \mathbb{Q}(\zeta_8, \sqrt{1+i}\,)$ generated by the square roots of the $\pi_j$ and $\sigma_j$.

For a proof of Proposition 7 we have to show that most of the torsors

$$\mathcal{T}^{(\phi)}(b_1): N^2 = b_1 M^4 + 4kM^2 e^2 + b_2 e^4, \quad b_1 b_2 = -4k^2, \tag{6}$$

do not have solutions[3]. To this end, we give necessary conditions for rational solvability that are stronger than solvability in all completions of $\mathbb{Q}$.

Write $k = k_1 k_2$ and consider the two possibilities $b_1 = k_1$ and $b_1 = 2k_1$. The next two lemmas will give necessary conditions for the solvability of (6).

**Lemma 10.** *Let $k = k_1 k_2$ and $b_1 = k_1 = p_1 \cdots p_r$, and write $k_1 = \kappa_1 \overline{\kappa}_1$ for some $\kappa_1 \in \mathbb{Z}[\sqrt{2}\,]$ with $\kappa_1 \equiv 1 \bmod 2$. If (6) has a solution, then*

- $b_1 \equiv 1 \bmod 16$;
- $[\kappa_1/\pi] = 1$ for all $\pi \mid k_2$;

---

[3] From now on, we only say solution when we actually mean non-trivial primitive integral solution.

- $[\kappa_1'/\pi] = 1$ *for all* $\pi \mid \kappa_1$, *where* $\kappa = \pi\kappa_1'$.

*Proof.* First observe that there exists some $\kappa_1$ with the desired properties since $k_1$ is a product of primes $\equiv 1 \bmod 8$. Next $\kappa_1$ is determined up to squares of units and conjugacy, which together with Lemma 1 implies that the symbols $[\kappa_1/\pi]$ do not depend on the choice of $\kappa_1$.

Next we first notice that $M$ must be odd: if it were even, $e$ would be odd, and this would imply $k_1 \equiv 3 \bmod 4$. From $(M, e) = 1$ we deduce that $(M, k_1) = 1$; we write $(M, k) = (M, k_2) = l$ for some integer $l > 0$ and put $M = lm$, $k_2 = lk_3$ and $N = k_1 ln$; this gives

$$k_1 n^2 = l^2 m^4 + 4k_2 m^2 e^2 - 4k_3^2 e^4 = (lm^2 + 2\varepsilon k_3 e^2)(lm^2 + 2\overline{\varepsilon} k_3 e^2),$$

where $\varepsilon = 1 + \sqrt{2}$ is the fundamental unit in $\mathbb{Z}[\sqrt{2}]$. Since any common divisor of the two factors divides both their sum and their difference, we see that they are coprime in $\mathbb{Z}[\sqrt{2}]$. Thus we find

$$lm^2 + 2\varepsilon k_3 e^2 = \kappa_1 \nu^2, \tag{7}$$

where $\kappa_1 \overline{\kappa}_1 = k_1$ and $\nu\overline{\nu} = n$. Reducing this equation modulo some $\pi \mid k_2$ gives $[\kappa_1/\pi] = 1$ (in fact, for $\pi \mid l$ we find $[\kappa_1/\pi] = [2\varepsilon k_3/\pi]$; but $2 = \sqrt{2}^2$ is a square, $[\varepsilon/\pi] = (1 + \sqrt{2}/p) = (1 + i/p) = 1$ by Lemma 5 and our assumptions, and $[k_3/\pi] = (k_3/p) = +1$ again by our assumptions. The case $\pi \mid k_3$ is similar), which proves our second claim.

Now we know that $l \equiv m^2 \equiv k_3 \equiv 1 \bmod 8$ in (7), which gives us either $\kappa_1 \nu^2 \equiv 1 \bmod 8$ or $\kappa_1 \nu^2 \equiv 1 + 2\varepsilon = \varepsilon^2 \bmod 8$ according as $e$ is even od odd. Thus $\kappa_1$ is congruent to a square modulo 8, that is, $\kappa_1 \equiv 1, 1 + 4\sqrt{2}, 3 \pm 2\sqrt{2} \bmod 8$ (in particular, $\kappa_1 \equiv 1 \bmod 2$). From this we deduce that $k_1 = N\kappa_1 \equiv 1 \bmod 16$ as claimed.

Subtracting (7) from its conjugate we get

$$4\sqrt{2} k_3 e^2 = \kappa_1 \nu^2 - \overline{\kappa}_1 \overline{\nu}^2. \tag{8}$$

Reducing (8) modulo some $\overline{\pi} \mid \overline{\kappa}_1$ we get $(2/p)_4 = [\sqrt{2}/\overline{\pi}] = [\kappa_1/\overline{\pi}] = [\kappa_1'/\overline{\pi}][\pi/\overline{\pi}] = [\kappa_1'/\overline{\pi}](2/p)_4$, where $\kappa_1 = \pi\kappa_1'$. This implies $[\kappa_1'/\overline{\pi}] = 1$. By Lemma 1, $[\kappa_1'/\overline{\pi}] = [\kappa_1'/\pi]$, and our proof is complete. $\square$

**Lemma 11.** *Let* $k = k_1 k_2$ *with* $k_1 = p_1 \cdots p_r$ *and* $b_1 = 2k_1$. *Write* $k_1 = \kappa_1 \overline{\kappa}_1$ *for some* $\kappa_1 \in \mathbb{Z}[\sqrt{2}]$ *with* $\kappa_1 \equiv 1 \bmod 2$. *If (6) has a solution, then*

- $[\kappa_1/\pi] = (2/p)_4$ *for all* $\pi \mid k_2$ *with* $p = N\pi$;

- $[\overline{\kappa}_1/\pi] = 1$ *for* $\pi \mid \kappa_1$.

*Proof.* Put $(M, k) = (M, k_2) = l$; then $N = 2k_1 ln$, $M = lm$, $k_2 = lk_3$ and

$$2k_1 n^2 = l^2 m^4 + 2k_2 m^2 e^2 - k_3^2 e^4 = (lm^2 + \varepsilon k_3 e^2)(lm^2 + \varepsilon' k_3 e^2).$$

The gcd of the two factors equals $\sqrt{2}$ (this follows from $l \equiv m \equiv e \equiv k_3 \equiv 1 \bmod 2$, since these congruences imply that both factors are $\equiv 1 + \varepsilon = 2 + \sqrt{2} \bmod 2$), hence we get $lm^2 + \varepsilon k_3 e^2 = \sqrt{2}\kappa_1 \nu^2$ and $lm^2 + \overline{\varepsilon} k_3 e^2 = -\sqrt{2}\overline{\kappa_1}\overline{\nu}^2$ for some $\kappa_1 \in \mathbb{Z}[\sqrt{2}]$ with $N\kappa_1 = k_1$. Reducing modulo some $\pi \mid p$ such that $p \mid k_2$ gives $(2/p)_4 [\kappa_1/\pi] = 1$.

Forming the difference as above shows that $2k_3 e^2 = \kappa_1 \nu^2 + \overline{\kappa_1}\overline{\nu}^2$; thus $[\kappa_1/\overline{\pi}] = [\overline{\kappa_1}/\pi] = 1$ for $\pi \mid \kappa_1$. □

**Some Special Cases.** Let us now look at the special case where $k = pq$, where $p \equiv q \equiv 1 \bmod 8$ are primes such that $(\frac{1+i}{p}) = (\frac{1+i}{q}) = +1$ and $(\frac{p}{q}) = +1$. Write $N\sigma = p$, $N\tau = q$ for some $\sigma, \tau \in \mathbb{Z}[\sqrt{2}]$ with $\sigma \equiv \tau \equiv 1 \bmod 2$ and $\chi(n) = (-1)^{(n-1)/8}$. We will examine the following eight cases:

| case | $\chi(p)$ | $\chi(q)$ | $[\frac{\sigma}{\tau}]$ | case | $\chi(p)$ | $\chi(q)$ | $[\frac{\sigma}{\tau}]$ |
|------|-----------|-----------|-------------------------|------|-----------|-----------|-------------------------|
| a)   | $-1$      | $-1$      | $-1$                    | e)   | $+1$      | $-1$      | $-1$                    |
| b)   | $-1$      | $-1$      | $+1$                    | f)   | $+1$      | $-1$      | $+1$                    |
| c)   | $-1$      | $+1$      | $-1$                    | g)   | $+1$      | $+1$      | $-1$                    |
| d)   | $-1$      | $+1$      | $+1$                    | h)   | $+1$      | $+1$      | $+1$                    |

In the last case h) we cannot draw any conclusions about the order of $\mathrm{III}(E_k/\mathbb{Q})[\phi]$, since the conditions of Lemmas 10 and 11 are all satisfied. In all other cases, however, it turns out that only one equation can have solutions: see the following table, which displays which of the equations (6) possibly have solutions according to the values of $\chi(p)$, $\chi(q)$ and $[\sigma/\tau]$.

| $b_1$ | conditions | a) | b) | c) | d) | e) | f) | g) |
|-------|------------|----|----|----|----|----|----|-----|
| $\pm 2$ | $\chi(p) = \chi(q) = 1$ | no | no | no | no | no | no | ? |
| $\pm p$ | $[\frac{\sigma}{\tau}] = +1, \chi(p) = 1$ | no | no | no | no | no | ? | no |
| $\pm 2p$ | $[\frac{\sigma}{\tau}] = \chi(q), \chi(p) = 1$ | no | no | no | no | ? | no | no |
| $\pm q$ | $[\frac{\sigma}{\tau}] = +1, \chi(q) = 1$ | no | no | no | ? | no | no | no |
| $\pm 2q$ | $[\frac{\sigma}{\tau}] = \chi(p), \chi(q) = 1$ | no | no | ? | no | no | no | no |
| $\pm pq$ | $\chi(pq) = 1, (\frac{\sigma}{\tau}) = +1$ | no | ? | no | no | no | no | no |
| $\pm 2pq$ | $[\frac{\sigma}{\tau}] = \chi(p) = \chi(q)$ | ? | no | no | no | no | no | no |

Here any entry "no" indicates that the corresponding torsor does not have a rational point. Observe that $\mathrm{III}(E/\mathbb{Q})[\phi]$ is generated by the classes of $2\mathbb{Q}^{\times 2}$ and $q\mathbb{Q}^{\times 2}$ whenever $\chi(q) = -1$, assuming that the equations marked with a "?" do have solutions in $\mathbb{Q}$ as predicted by the parity conjecture. Table 1 gives a few numerical solutions of the torsors $\mathcal{T}^{(\phi)}(b_1) : b_1 n^2 = m^4 + 2c_1 m^2 e^2 - c_1^2 e^4$, where $2k = b_1 c_1$.

The solutions for the torsors $\mathcal{T}^{(\phi)}(2)$ and $\mathcal{T}^{(\phi)}(p)$ in case h) were not found by a direct (and very naive) search: rather, we used the smaller solutions found for $b_1 = q, pq$ and $2pq$ to compute them from the group law on torsors that we have explained in Section 1.

Table 1:

| case | $p$ | $q$ | $b_1$ | $n$ | $m$ | $e$ | $x$ | $y$ |
|---|---|---|---|---|---|---|---|---|
| a) | 41 | 409 | $2pq$ | 31 | 71 | 27 | $\frac{169065058}{729}$ | $\frac{2475679174244}{19683}$ |
| b) | 41 | 569 | $pq$ | 7 | 31 | 30 | $\frac{22419169}{900}$ | $\frac{118100566297}{27000}$ |
| c) | 41 | 353 | $2q$ | 1927 | 205 | 17 | $\frac{29669650}{289}$ | $\frac{196899665260}{4913}$ |
| d) | 457 | 593 | $q$ | 99169 | 1371 | 29 | $\frac{1114627113}{841}$ | $\frac{47810443842651}{24389}$ |
| e) | 353 | 41 | $2p$ | 1927 | 205 | 17 | $\frac{29669650}{289}$ | $\frac{196899665260}{4913}$ |
| f) | 113 | 41 | $p$ | 41 | 41 | 7 | $\frac{189953}{49}$ | $\frac{21464689}{343}$ |
| g) | 113 | 257 | 2 | 100487 | 771 | 7 | $\frac{1188882}{49}$ | $\frac{309901908}{343}$ |
| h) | 113 | 337 | 2 | 4971761 | 3033 | 23 | $\frac{18398178}{529}$ | $\frac{60317404452}{12167}$ |
|  |  |  | $p$ | 167279 | 4047 | 242 | $\frac{1850737617}{58564}$ | $\frac{8644333524897}{14172488}$ |
|  |  |  | $q$ | 791 | 113 | 3 | $\frac{4303153}{9}$ | $\frac{10151137927}{27}$ |

Now assume that $\chi(p) = \chi(q) = \chi(r) = -1$ and consider the following cases:

| case | $[\pi/\lambda]$ | $[\pi/\rho]$ | $[\rho/\lambda]$ | case | $[\pi/\lambda]$ | $[\pi/\rho]$ | $[\rho/\lambda]$ |
|---|---|---|---|---|---|---|---|
| a) | $-1$ | $-1$ | $-1$ | e) | $+1$ | $-1$ | $-1$ |
| b) | $-1$ | $-1$ | $+1$ | f) | $+1$ | $-1$ | $+1$ |
| c) | $-1$ | $+1$ | $-1$ | g) | $+1$ | $+1$ | $-1$ |
| d) | $-1$ | $+1$ | $+1$ | h) | $+1$ | $+1$ | $+1$ |

Then it turns out that all of the torsors $\mathcal{T}^{(\phi)}(b)$ with $b \neq \pm 1$ are non-trivial in the cases a), d), f) and g), that is, in these cases (note that they are characterized by the condition $[\pi/\lambda][\lambda/\rho][\rho/\pi] = -1$) we have $\#\mathrm{III}(E/\mathbb{Q})[\phi] = \frac{1}{2}\# S^{(\phi)}(E/\mathbb{Q})$. In the other cases we find results that are similar to the case $k = pq$; see Table 2.

**Proof of Proposition 7.** Our first observation is that it suffices to consider positive values of $b_1$: since $-1\mathbb{Q}^{\times 2} \in W(E/\mathbb{Q})$, the torsor $\mathcal{T}^{(\phi)}(b_1)$ possesses a solution if and only if $\mathcal{T}^{(\phi)}(-b_1)$ does.

Let us first deal with the case where $b_1 = k_1$ is a product of $r$ primes $p_j$. If $r$ is odd, then there exists a $\pi \mid k_2$ unless $r = t$, hence $[\kappa_1/\pi] = (-1)^r = -1$, and now Lemma 10 says that $\mathcal{T}^{(\phi)}(b_1)$ has no rational solution. If $r = t$, on the other hand, then $b_1 = k \equiv 9 \bmod 16$, and again Lemma 10 prohibits a rational solution.

Table 2:

| $b_1$ | conditions | b) | c) | e) | h) |
|---|---|---|---|---|---|
| $p$ | $[\pi/\lambda] = [\pi/\rho] = \chi(p) = 1$ | no | no | no | no |
| $q$ | $[\lambda/\pi] = [\lambda/\rho] = \chi(q) = 1$ | no | no | no | no |
| $r$ | $[\rho/\pi] = [\rho/\lambda] = \chi(r) = 1$ | no | no | no | no |
| $pq$ | $[\pi\lambda/\rho] = [\lambda/\pi] = \chi(pq) = 1$ | no | no | ? | ? |
| $qr$ | $[\lambda\rho/\pi] = [\lambda/\rho] = \chi(qr) = 1$ | ? | no | no | ? |
| $rp$ | $[\pi\rho/\lambda] = [\pi/\rho] = \chi(rp) = 1$ | no | ? | no | ? |
| $pqr$ | $[\lambda\rho/\pi] = [\pi\rho/\lambda] = \chi(pqr) = 1$ | no | no | no | no |
| $2$ | $\chi(p) = \chi(q) = \chi(r) = 1$ | no | no | no | no |
| $2p$ | $\chi(p) = 1, \chi(q) = [\lambda/\pi], \chi(r) = [\pi/\rho]$ | no | no | no | no |
| $2q$ | $\chi(q) = 1, \chi(r) = [\lambda/\rho], \chi(p) = [\lambda/\pi]$ | no | no | no | no |
| $2r$ | $\chi(r) = 1, \chi(p) = [\pi/\rho], \chi(q) = [\lambda/\rho]$ | no | no | no | no |
| $2pq$ | $\chi(r) = [\pi\lambda/\rho], \chi(p) = \chi(q) = [\lambda/\pi]$ | ? | ? | no | no |
| $2qr$ | $\chi(p) = [\lambda\rho/\pi], \chi(q) = \chi(r) = [\lambda/\rho]$ | no | ? | ? | no |
| $2rp$ | $\chi(q) = [\pi\rho/\lambda], \chi(p) = \chi(r) = [\pi/\rho]$ | ? | no | ? | no |
| $2pqr$ | $\chi(p) = [\lambda\rho/\pi], \chi(q) = [\pi\rho/\lambda], \chi(r) = [\pi\lambda/\rho]$ | no | no | no | no |

If $r$ is even and $\neq 0$, then there is some $\pi \mid \kappa_1$, and we have $[\overline{\kappa}'_1/\pi] = (-1)^{r-1} = -1$. Invoking Lemma 10 again, we see that rational solvability implies that $r = 0$, that is, $b_1 = \pm 1$.

Next consider the case $b_1 = 2k_1$ with $k_1$ as before. If $r$ is even and $\neq t$, then there is a prime $\pi$ with norm $p \mid k_2$, and we have $\chi(p) = -1$ and $[\kappa_1/\pi] = (-1)^r = +1$. Thus $\mathcal{T}^{(\phi)}(b_1)$ does not have a rational solution for even $r$ by Lemma 11 unless $t$ is even and $b_1 = \pm 2k$. If $r$ is odd, let $\pi$ be a prime dividing $\kappa_1$. Then solvability of our torsor implies that $1 = [\overline{\kappa}_1/\pi] = [\overline{\pi}/\pi][\overline{\kappa}'_1/\pi] = \chi(p)(-1)^{r-1} = (-1)^r$: contradiction.

**The Isogenous Curve $\widehat{E}_k$.** Next we study the Tate-Shafarevich group of the isogenous curve $\widehat{E}$. Here, our torsor

$$\mathcal{T}^{(\psi)}(b_1) : N^2 = b_1 M^4 - 2kM^2 e^2 + b_2 e^4, \quad b_1 b_2 = 2k^2, \tag{9}$$

factors over the Gaussian integers $\mathbb{Z}[i]$. Since $2\mathbb{Q}^{\times 2} \in W(\widehat{E}/\mathbb{Q})$, it is sufficient to consider odd values of $b_1$.

So let us take $b_1 = k_1$, with $k$, $k_1$, etc. as in the discussion above. We immediately find that $M$ must be odd. Write $(M,k) = l$. Then $M = lm$, $n = k_1 ln$ and $k_2 = lk_3$, and we find

$$k_1 n^2 = l^2 m^4 - 2k_2 m^2 e^2 + 2k_3^2 e^4 = (lm^2 - (1+i)k_3 e^2)(lm^2 - (1-i)k_3 e^2).$$

Since $l$ and $m$ are odd, the factors on the right hand side are coprime, and unique factorization in $\mathbb{Z}[i]$ gives

$$\kappa_1 \nu^2 = lm^2 - (1+i)k_3 e^2, \tag{10}$$

where $\kappa_1 \in \mathbb{Z}[i]$ has norm $k_1$. If $\sigma \in \mathbb{Z}[i]$ is a prime with norm $p \mid k_2$, then equation (10) implies $[\kappa_1/\sigma] = +1$ (note that $(\frac{1+i}{p}) = +1$ by assumption, and that $[\kappa_1/\sigma]$ does not depend on the choice of $\kappa_1$ since $[i/\sigma] = +1$ for any prime $\sigma$ with norm $\equiv 1 \bmod 8$).

Subtracting (10) from its conjugate and reducing modulo primes $\sigma \mid \overline{\kappa}_1$ gives $[\overline{\kappa}_1/\sigma] = +1$. We have proved

**Lemma 12.** *If the torsor (9) has a solution then*

- $[\kappa_1/\sigma] = +1$ *for all primes* $\sigma \in \mathbb{Z}[i]$ *with norm dividing $k_2$;*
- $[\kappa_1/\sigma] = +1$ *for all primes* $\sigma \mid \overline{\kappa}_1$,

*where $\kappa_1 \overline{\kappa}_1 = k_1$ for some $\kappa_1 \in \mathbb{Z}[i]$.*

Now we can prove Proposition 8: assume that $b_1 = k_1$ is the product of $r$ primes $p_j$. If $r$ is odd and $k_2 \neq 1$, then $[\kappa_1/\sigma] = (-1)^r = -1$, and $\mathcal{T}^{(\psi)}(b_1)$ is not solvable. In other words: if $r$ is odd and $\mathcal{T}^{(\psi)}(b_1)$ has a solution then $r = t$ and $b_1 = k_1 = k$. If $r \neq 0$ is even, then there is a $\sigma \mid \overline{\kappa}_1$, and we have $[\kappa_1/\sigma] = [\overline{\sigma}/\sigma](-1)^{r-1} = -1$ since $[\overline{\sigma}/\sigma] = +1$. Thus $\mathcal{T}^{(\psi)}(b_1)$ is not solvable.

As a corollary we note that there exist elliptic curves $E_k$ that are 2-isogenous to $\widehat{E}_k$ such that the 2-ranks of $\Sha(E_k/\mathbb{Q})$ and $\Sha(\widehat{E}_k/\mathbb{Q})$ both grow arbitrarily large. This is a consequence of the embedding $\Sha(E_k/\mathbb{Q})[\phi] \longrightarrow \Sha(E_k/\mathbb{Q})$ discussed in Section 5.

## 4. More Examples

In this section we apply the technique explained above to the curves $E_{-k}$ with $k$ as in Proposition 6 without giving the technical details. Under the assumptions of Proposition 6, the Selmer groups for the $E_{-k}$ are

$$\begin{aligned} S^{(\phi)}(E_{-k}/\mathbb{Q}) &= \{b_1 \mathbb{Q}^{\times 2} : b_1 \text{ squarefree}, b_1 \mid k\}, \\ S^{(\psi)}(\widehat{E}_{-k}/\mathbb{Q}) &= \{b_1 \mathbb{Q}^{\times 2} : b_1 \text{ squarefree}, 0 < b_1 \mid 2k\}. \end{aligned}$$

For the calculation of the Tate-Shafarevich groups we need the analogs of our lemmas; Lemma 10 holds as stated except that we can no longer conclude that $k_1 \equiv 1 \bmod 16$: in fact, we only get that $\kappa_1 \nu^2 \equiv 1, 1 - 2\varepsilon \bmod 8$ (as opposed to $\kappa_1 \nu^2 \equiv 1, 1 + 2\varepsilon \bmod 8$ for positive $k$), and $1 - 2\varepsilon = -1 - 2\sqrt{2}$ is not congruent to $\pm \xi^2 \bmod 8$ for $\xi \in \mathbb{Z}[\sqrt{2}]$. In fact it can be shown that the conditions we get for $\kappa_1$ are equivalent to $[1 + \sqrt{2}/\kappa_1] = 1$ (see Lemma 5), that is, they are satisfied automatically since we have to assume that $(1 + \sqrt{2}/p_j) = (1 + i/p_j) = 1$ in order to make our torsors everywhere locally solvable (see Proposition 6).

On the other hand, we find that Lemma 11 holds, and that in addition we have $k_1 \equiv 1 \bmod 16$. Finally, Lemma 12 goes through unchanged. Using these lemmas, the analogs of our results in Section 3 are easily derived:

**Proposition 13.** *Consider the curve $E = E_{-k}$ where $k = p_1 \cdots p_t$ is a product of primes $p_i \equiv 1 \bmod 8$, and assume that $k \equiv 9 \bmod 16$ if $t$ is odd. Suppose moreover that $(p_i/p_j) = +1$ for $i \neq j$ and $(\frac{1+i}{p}) = +1$ for all $p \mid k$; write $p_i = \pi_i \overline{\pi}_i$ for elements $\pi_i, \overline{\pi}_i \in \mathbb{Z}[\sqrt{2}]$ with $\pi_i \equiv 1 \bmod 2$, and assume in addition that $[\pi_i/\pi_j] = -1$ for $i \neq j$. Then $W(E/\mathbb{Q}) = \{\mathbb{Q}^{\times 2}, -\mathbb{Q}^{\times 2}\}$ if $t$ is even, and $W(E/\mathbb{Q}) = \{\mathbb{Q}^{\times 2}, -\mathbb{Q}^{\times 2}\}$ or $W(E/\mathbb{Q}) = \{\pm \mathbb{Q}^{\times 2}, \pm k \mathbb{Q}^{\times 2}\}$ if $t$ is odd.*

*In particular, we have $\#\text{III}(E/\mathbb{Q})[\phi] = \frac{1}{2} \#S^{(\phi)}(E/\mathbb{Q}) = 2^t$ if $t$ is even, and $\#\text{III}(E/\mathbb{Q})[\phi] \geq \frac{1}{4} \#S^{(\phi)}(E/\mathbb{Q}) = 2^{t-1}$ if $t$ is odd.*

Here is the corresponding result for $\text{III}(\widehat{E}/\mathbb{Q})[\psi]$:

**Proposition 14.** *Consider the elliptic curve $E = E_{-k}$ where $k = p_1 \cdots p_t$ is a product of primes $p_i \equiv 1 \bmod 8$ with $(p_i/p_j) = +1$ for $i \neq j$ and $(\frac{1+i}{p}) = +1$ for all $p \mid k$; write $p_i = \sigma_i \overline{\sigma}_i$ for elements $\sigma_i, \overline{\sigma}_i \in \mathbb{Z}[i]$, and assume in addition that $[\sigma_i/\sigma_j] = -1$ for $i \neq j$. Then $W(\widehat{E}/\mathbb{Q}) = \{\mathbb{Q}^{\times 2}, 2\mathbb{Q}^{\times 2}\}$ if $t$ is even, and $W(\widehat{E}/\mathbb{Q}) \subseteq \{\mathbb{Q}^{\times 2}, 2\mathbb{Q}^{\times 2}, k\mathbb{Q}^{\times 2}, 2k\mathbb{Q}^{\times 2}\}$ if $t$ is odd.*

*In particular, we have $\#\text{III}(\widehat{E}/\mathbb{Q})[\psi] \geq \frac{1}{4} \#S^{(\psi)}(\widehat{E}/\mathbb{Q}) = 2^{t-1}$ if $t$ is odd, and $\#\text{III}(\widehat{E}/\mathbb{Q})[\psi] = \frac{1}{2} \#S^{(\psi)}(\widehat{E}/\mathbb{Q}) = 2^t$ if $t$ is even.*

Again we can combine these results and get

**Theorem 15.** *If the conditions in Propositions 7 and 8 are satisfied, then $E_{-k}(\mathbb{Q})$ has rank $\leq 2$ if $t$ is odd, and rank $0$ if $t$ is even. If the parity conjecture holds and $t$ is odd, then either $\operatorname{rank} E_{-k}(\mathbb{Q}) = 2$ and the inequalities in Propositions 13 and 14 are equalities, or $\operatorname{rank} E_{-k}(\mathbb{Q}) = 0$, and then $W(E_{-k}/\mathbb{Q}) = \{\mathbb{Q}^{\times 2}, -\mathbb{Q}^{\times 2}\}$ and $W(\widehat{E}_{-k}/\mathbb{Q}) = \{\mathbb{Q}^{\times 2}, 2\mathbb{Q}^{\times 2}\}$.*

**Remark.** As a conclusion we remark that there is still a lot of work to do: although we computed the 2-ranks of $\text{III}(E/\mathbb{Q})[\phi]$ and $\text{III}(\widehat{E}/\mathbb{Q})[\psi]$ for some $E = E_k$ and $E = E_{-k}$, we do not know yet what the rank of $\text{III}(E/\mathbb{Q})[2]$ or $\text{III}(\widehat{E}/\mathbb{Q})[2]$ is. Techniques explained by Birch and Swinnerton-Dyer [2], McGuinness [12] and Cremona [6] should suffice for computing the exact 2-rank of $\text{III}(E/\mathbb{Q})$, and using Cassels [5] it should be possible to compute the 4-rank as well, at least for explicitly given curves.

**Remark.** We also note that our results can of course be used to search for curves $E_k$ with large rank: one only need consider products $k = p_1 \cdots p_r \equiv 1 \bmod 16$ of primes $p_i \equiv 1 \bmod 8$ such that the quadratic residue symbols in Propositions 7 and 8 are $+1$ instead of $-1$.

## 5. A Conjecture

The computations made in this paper suggest a conjecture that we are going to discuss now.

**Conjecture A.** Let $E/\mathbb{Q}$ be a modular elliptic curve with a rational 2-torsion point $\neq \mathcal{O}$. Then

$$\text{rank } S^{(\phi)}(E/\mathbb{Q}) + \text{rank } S^{(\psi)}(\widehat{E}/\mathbb{Q}) \equiv v \bmod 2,$$

where $w = (-1)^v$ is the root number, that is, the sign in the functional equation of the complete Hasse-Weil $L$-function $\xi(E, s)$ of $E$.

This conjecture is nowhere near as deep as e.g. the parity conjecture: the latter predicts the existence of rational points on certain elliptic curves, whereas Conjecture A only compares local data with the sign of the functional equation. In fact, since the Selmer groups for the curves $E_k$ are explicitly given in [13], it appears to be only a combinatorial problem to verify this conjecture for these curves; apart from the curves $E_{\pm k}$ treated in Proposition 6 and below I checked it for all $E_k$ with $k = u$ and $k = up$, where $u \in \{\pm 1, \pm 2\}$ and where $p$ is any prime. The next thing to do is to prove that Conjecture A is compatible with quadratic twists (note we know how the right hand side behaves under quadratic twists).

Conjecture A predicts that rank $\text{III}(E/\mathbb{Q})[\phi]$ and rank $\text{III}(\widehat{E}/\mathbb{Q})[\psi]$ are either both odd or both even, assuming the parity conjecture. Actually, Conjecture A and the parity conjecture are equivalent under the assumption that $\#\text{III}(E/\mathbb{Q})[2]$ is a square. For a proof, consider the commutative exact diagram (see e.g. Cremona [6])

$$\begin{array}{ccccccccc}
K & \longrightarrow & W(E/\mathbb{Q}) & \longrightarrow & E(\mathbb{Q})/2E(\mathbb{Q}) & \longrightarrow & W(\widehat{E}/\mathbb{Q}) & \longrightarrow & 0 \\
\downarrow & & \downarrow & & \downarrow & & \downarrow & & \downarrow \\
K & \longrightarrow & S^{(\phi)}(E/\mathbb{Q}) & \longrightarrow & S^{(2)}(E/\mathbb{Q}) & \longrightarrow & S^{(\psi)}(\widehat{E}/\mathbb{Q}) & \longrightarrow & \widehat{C} \\
\downarrow & & \downarrow & & \downarrow & & \downarrow & & \downarrow \\
0 & \longrightarrow & \text{III}(E/\mathbb{Q})[\phi] & \longrightarrow & \text{III}(E/\mathbb{Q})[2] & \longrightarrow & \text{III}(\widehat{E}/\mathbb{Q})[\psi] & \longrightarrow & \widehat{C}
\end{array}$$

Here

$$K = \ker \, [W(E/\mathbb{Q}) \longrightarrow E(\mathbb{Q})/2E(\mathbb{Q})] \simeq \ker \, [S^{(\phi)}(E/\mathbb{Q}) \longrightarrow S^{(2)}(E/\mathbb{Q})]$$

is a group of order $4/\#E(\mathbb{Q})[2]$, and Cassels proved that

$$\widehat{C} = \text{cok} \, [S^{(2)}(E/\mathbb{Q}) \longrightarrow S^{(\psi)}(\widehat{E}/\mathbb{Q})] \simeq \text{cok} \, [\text{III}(E/\mathbb{Q})[2] \longrightarrow \text{III}(\widehat{E}/\mathbb{Q})[\psi]]$$

has square order. Thus $\#S^{(\phi)}(E/\mathbb{Q}) \cdot \#S^{(\psi)}(\widehat{E}/\mathbb{Q}) = \#K \cdot \#S^{(2)}(E/\mathbb{Q}) \cdot \#\widehat{C}$, and since $\#S^{(2)}(E/\mathbb{Q}) = 2^r \#E(\mathbb{Q})[2]/\#\text{III}(E/\mathbb{Q})[2]$ this gives

$$\#S^{(\phi)}(E/\mathbb{Q})\#S^{(\psi)}(\widehat{E}/\mathbb{Q}) = 2^r \#\widehat{C}/\#\text{III}(E/\mathbb{Q})[2],$$

that is,

$$\text{rank } S^{(\phi)}(E/\mathbb{Q}) + \text{rank } S^{(\psi)}(\widehat{E}/\mathbb{Q}) \equiv r \bmod 2$$

if $\#\text{III}(E/\mathbb{Q})[2]$ is a square, and taking this for granted, Conjecture A is equivalent to the parity conjecture as claimed. As I recently discovered, Conjecture A was proved by Birch & Stephens [1] for the elliptic curves of type $y^2 = x^3 + Dx$ and $y^2 = x^3 + D$.

The diagram plus the fact that $\#C$ is a square also imply that the product $\#\text{III}(E/\mathbb{Q})[\phi] \cdot \#\text{III}(\widehat{E}/\mathbb{Q})[\psi]$ is a square if and only $\#\text{III}(E/\mathbb{Q})[2]$ is one. Since our calculations produced Tate-Shafarevich groups whose orders are squares one is tempted to conjecture that this is always true. But as Oisin McGuinness explained to me, it does not follow from the existence of the Cassels pairing on $\text{III}(E/\mathbb{Q})$ that the orders of Tate-Shafarevich groups corresponding to isogenies of non-square degree are perfect squares. The first example for such a behaviour I found were the groups $\text{III}(E/k)[\phi]$ and $\text{III}(\widehat{E}/k)[\psi]$ for the curve $E : y^2 = x^3 - 68x$ over the quadratic number field $k = \mathbb{Q}(\sqrt{-43})$, both of order 2 (see [10]). More recently, I have found infinitely many examples over $\mathbb{Q}$ using elliptic curves with three rational points of order 2 (see [11]).

# Acknowledgements

I thank J.W.S. Cassels and W. Stein for pointing out the references to the papers [3], [4] and [8], and O. McGuinness for sending me an extract of [12] as well as for his explanation why $\#\text{III}(E/\mathbb{Q})[\phi]$ cannot be expected to be a square for 2-isogenies $\phi$. Last not least I gratefully acknowledge the financial support of the DFG.